%% file: top.tex
\documentclass[11pt]{amsart}

\usepackage{xparse}
\usepackage{etoolbox}
\usepackage{aliascnt}

\RequirePackage{doi}

\usepackage{courier}
\usepackage[T1]{fontenc}

\usepackage[scaled]{helvet}
\usepackage[mathbf]{euler}
\usepackage[driver=pdftex,margin=3cm,heightrounded=true,centering]{geometry}

\usepackage{enumerate}
\usepackage{amssymb}
\usepackage{amsopn}
\usepackage{amsmath}
\usepackage{mathtools}

\usepackage{tikz}
\usepackage{tikz-cd}
\usepackage{ifthen}
\usepackage{graphics}
\usepackage{hyperref}
\hypersetup{
	colorlinks=true,
	linkcolor=blue,
	filecolor=blue,
	citecolor = blue,
	urlcolor=blue,
	bookmarksdepth=3,
	breaklinks=true,
}

\usepackage{amsthm}
\usepackage{topthm}

\usepackage[textsize = footnotesize]{todonotes}
\usepackage{comment}
\author{Thorben Kastenholz}
\thanks{}
\date{\today}
\title{Simplcial volume of manifolds fibering with connected structure group}
\address{Karlsruher Institut f\"ur Technologie, Englerstraße 2, 76131
Karlsruhe, Germany}
\email{thorben.kastenholz@kit.edu}
\begin{document}
\input{commands.tex}

\input{Abstract.tex}
\maketitle
\section{Introduction}
\input{Section/Introduction.tex}
\section{Fiber bundles and structure groups}
\input{Section/FBandSG.tex}
\section{Proof of \autoref{thm:RationalProductVanishing}}
\input{Section/Proof.tex}
\section{Finding small representatives}
\input{Section/Representatives.tex}

\bibliography{sources}
\bibliographystyle{plain}
\end{document}

%% file: commands.tex
\newcommand{\introduce}[1]
  {\textbf{#1}}
\newcommand{\tk}[1]{\todo[size=\tiny,color=green!40]{TK: #1}}
\newcommand\blfootnote[1]{%
  \begingroup
  \renewcommand\thefootnote{}\footnote{#1}%
  \addtocounter{footnote}{-1}%
  \endgroup
}

\newcommand{\apply}[2]
  {{#1}\!\left({#2}\right)}
\newcommand{\at}[2]
  {\left.{#1}\right\rvert_{#2}}
\newcommand{\Identity}%
  {\mathrm{Id}}
\newcommand{\NaturalNumbers}%
  {\mathbf{N}}
\newcommand{\Integers}%
  {\mathbf{Z}}
\newcommand{\Rationals}%
  {\mathbf{Q}}
\newcommand{\Reals}%
  {\mathbf{R}}
  \newcommand{\ComplexNumbers}%
  {\mathbf{C}}
\newcommand{\AbstractProjection}[1] 
  {p_{#1}}
\newcommand{\RealPart}[1]
  {\apply{\operatorname{Re}}{#1}}
\newcommand{\ImaginaryPart}[1]
  {\apply{\operatorname{Im}}{#1}}
\newcommand{\Floor}[1]
  {\left \lfloor #1 \right \rfloor}
\newcommand{\Norm}[1]
  {\left|\left|#1\right|\right|}
\newcommand{\AbsoluteValue}[1]
  {\left|#1\right|}
\newcommand{\GeneralLinearGroup}[1]
  {\apply{\mathop{GL}_{#1}}{\Reals}}

\newcommand{\Surface}[1]
  {\Sigma_{#1}}
\newcommand{\SurfaceGroup}[1]
  {S_{#1}}
\newcommand{\Manifold}%
  {M}
\newcommand{\Tangenbundle}[1]
  {T#1}
\newcommand{\NormalBundle}[2]
  {N_{#2}#1}
\newcommand{\FiberTransferHomology}[1]
  {#1^{!}}
\newcommand{\Submanifold}%
  {S}
\newcommand{\TubularNeighborhood}[1]
  {U_{#1}}
\newcommand{\ManifoldAlternative}%
  {N}
 \newcommand{\ManifoldAuxiliary}%
  {K}
\newcommand{\NullBordism}%
  {W}
\newcommand{\Bordism}
  {P}
 \newcommand{\ManifoldFiber}%
 {F}
 \newcommand{\ManifoldTotal}%
 {E}
  \newcommand{\ManifoldBase}%
  {B}
\newcommand{\SmoothMap}%
  {\phi}
  \newcommand{\MorseFunction}%
  {f}
\newcommand{\Diffeomorphism}%
  {\Phi}
\newcommand{\Dimension}%
  {d}
  \newcommand{\HalfDimension}%
  {n}
\newcommand{\FundamentalClass}[1]
  {\left[#1\right]}
\newcommand{\Interval}%
  {I}
\newcommand{\Ball}[1]
  {D^{#1}}
\newcommand{\Sphere}[1] 
  {S^{#1}}
\newcommand{\Torus}[1]
  {T^{#1}}
\newcommand{\SimplicialVolume}[1]
  {\lvert \lvert #1 \rvert \rvert}
\newcommand{\ellone}%
  {\ell_{1}}
\newcommand{\Boundary}[1]
  {\partial #1}
\newcommand{\ComplexOfEmbeddings}[1]
  {\apply{K}{#1}}
\newcommand{\GenusOf}[1]
  {\apply{\Genus}{#1}}
\newcommand{\StableGenusOf}[1]
  {\apply{\overline{\Genus}}{#1}}
\newcommand{\Surgery}%
  {\natural}
\newcommand{\HandleEmbedding}%
  {\Phi}
\newcommand{\Singularity}%
  {X}
\newcommand{\FundamentalGroupShortHand}
  {\pi}

\newcommand{\Metric}%
  {g}
\newcommand{\IsometryGroup}[2]
  {\apply{\mathrm{Isom}}{#1,#2}}

\newcommand{\Diff}[1]
  {\mathrm{Diff}\!\left(#1\right)}
\newcommand{\DiffGroup}[1]
  {\mathrm{Diff}^{B\Group}\!\left(#1\right)}
\newcommand{\DiffZero}[1]
  {\mathrm{Diff}_0\!\left(#1\right)}
\newcommand{\DiffOne}[1]
  {\widetilde{\mathrm{Diff}}_0\!\left(#1\right)}
\newcommand{\HomeoGroup}[1]
  {\apply{\mathrm{Homeo}}{#1}}
\newcommand{\HomeoCompactGroup}[1]
  {\apply{\mathrm{Homeo}_{c}}{#1}}
\newcommand{\HomeoLowerGroup}[1]
  {\apply{\mathrm{Homeo}^{\geq}}{#1}}

\newcommand{\FiberingSpace}%
  {E}
\newcommand{\FiberingProjektion}[1] 
  {\pi_{#1}}
\newcommand{\Fiber}%
  {F}
\newcommand{\FiberDimension}%
  {f}
\newcommand{\Base}%
  {B}
\newcommand{\ClutchingFunction}[1] 
  {\varphi_{#1}}
\newcommand{\TautologicalBundle}[1]
  {\gamma_{#1}}
\newcommand{\Automorphisms}[1]
  {\apply{\mathop{Aut}}{#1}}
\newcommand{\AutomorphismsConnected}[1]
  {\apply{\mathop{Aut}_0}{#1}}

\newcommand{\Group}%
  {G}
\newcommand{\CompactGroup}
  {K}
\newcommand{\Subgroup}
  {H}
\newcommand{\AmenableGroup}
  {A}
\newcommand{\GroupElement}%
  {g}
\newcommand{\Genus}%
  {g}
\newcommand{\QuadraticModule}%
  {\mathbf{M}}
\newcommand{\WittIndex}[1]
  {\apply{\Genus}{#1}}
\newcommand{\StableWittIndex}[1]
  {\apply{\overline{\Genus}}{#1}}
\newcommand{\ComplexOfHyperbolicInclusions}[1]
  {\apply{K^{a}}{#1}}
\newcommand{\ChainContraction}[1]
  {H_{#1}}
\newcommand{\FreeGroup}[1]
  {F_{#1}}

\newcommand{\HomotopyFiber}[1]
  {\apply{\mathop{HoFib}}{#1}}

\newcommand{\HomologyClass}%
  {\alpha}
\newcommand{\FramedBordism}[2]
  {\apply{\Omega^{fr}_{#1}}{#2}}
\newcommand{\HomologyOfSpaceObject}[3]
  {\apply{H_{#1}}{#2 ; #3}}
\newcommand{\CohomologyOfSpaceObject}[3]
  {\apply{H^{#1}}{#2 ; #3}}
\newcommand{\BoundedCohomologyOfSpaceObject}[3]
  {\apply{H^{#1}_{\text{b}}}{#2 ; #3}}
\newcommand{\BoundedCohomologyOfSimplicialObject}[3]
  {\apply{H^{#1}_{\text{b, s}}}{#2 ; #3}}
\newcommand{\HomologyOfSpaceMorphism}[1]
  {{#1}_{\ast}}
\newcommand{\HomologyOfGroupObject}[3]
  {\apply{H_{#1}}{#2; #3}}
\newcommand{\HomologyOfGroupMorphism}[2]
  {{#1}_{\ast}}
\newcommand{\HomologyOfSpacePairObject}[3]
  {\apply{H_{#1}}{{#2},{#3}}}
\newcommand{\Multiple}%
  {\lambda}
\newcommand{\AbstractChainComplex}[1]
  {C_{#1}}
\newcommand{\SimplicialChainComplex}[2]
  {\apply{\AbstractChainComplex{#1}}{#2}}
\newcommand{\Chain}%
  {\sigma}
\newcommand{\Diameter}[1]
  {\apply{\mathrm{Diam}}{#1}}
\newcommand{\BoundingChain}%
  {\rho}
\newcommand{\BoundaryChainComplex}%
  {\partial}
\newcommand{\FiberTransfer}[1]
  {#1^{!}}
\newcommand{\TopologicalSpace}%
  {X}
\newcommand{\MappingTorus}[1]
  {T_{#1}}
\newcommand{\Inner}[1]
  {\mathring{#1}}
\newcommand{\Point}%
  {\ast}
\newcommand{\Loop}%
  {\gamma}
\newcommand{\ContinuousMap}%
  {f}
  \newcommand{\ContinuousMapALT}%
  {g}
 \newcommand{\maps}%
  {\ensuremath{\text{maps}}}
\newcommand{\HomotopyGroupOfObject}[3]
  {\apply{\pi_{#1}}{{#2},{#3}}}
\newcommand{\HomotopyGroupOfPairObject}[4]
  {\apply{\pi_{#1}}{{#2},{#3},{#4}}}
\newcommand{\HomotopyGroupMorphism}[1] 
  {{#1}_{\ast}}
\newcommand{\EMSpace}[2]
  {\apply{K}{{#1},{#2}}}
\newcommand{\ClassifyingSpace}[1] 
  {B#1}
\newcommand{\ClassifyingMap}[1]
  {\zeta_{#1}}
\newcommand{\UniversalCovering}[1] 
  {\widetilde{#1}}
\newcommand{\UniversalCoveringMap}[1] 
  {\widetilde{#1}}
\newcommand{\FundamentalCycle}[1]
  {\sigma_{#1}}
\newcommand{\CohomologyClass}%
  {\phi}

\newcommand{\SimplicialComplex}
  {X}
\newcommand{\AuxSimplicialComplex}
  {K}
\newcommand{\Subcomplex}
  {Y}
\newcommand{\AuxSubcomplex}
  {L}
\newcommand{\Simplex}[1]
  {\sigma_{#1}}
\newcommand{\Link}[2]
  {\apply{\text{Lk}_{#1}}{#2}}
\newcommand{\Star}[2]
  {\apply{\text{St}_{#1}}{#2}}
\newcommand{\BoundaryIndexSimplex}[2]
  {\apply{\partial_{#1}}{#2}}
\newcommand{\BoundarySimplex}%
  {\partial}
\newcommand{\StandardSimplex}[1]
  {\Delta_{#1}}
\newcommand{\vertex}
  {v}
\newcommand{\GeometricRealization}[1]
  {\left\lvert #1 \right\rvert}
\newcommand{\BoundHomotopy}
  {N}
\newcommand{\Horn}[2]
  {\Lambda^{#1}_{#2}}
\newcommand{\Triangulation}[1]
  {T_{#1}}

%% file: Abstract.tex
\begin{abstract}
  In this note we investigate the simplicial volume of fiber bundles with
  connected structure group. We are able to show that if the structure group is
  either compact or a Lie group, or if the fiber is aspherical that the
  simplicial volume of the total space agrees with the simplicial volume of the
  trivial bundle.
\end{abstract}

%% file: Section/Introduction.tex
For a homology class $\alpha$ of a topological space, one can define its
$\ellone$-norm via
\[
  \Norm{\alpha}
  =
  \inf
  \left\{
    \sum_{i}
    \AbsoluteValue{\lambda_i}
  \mid
  \alpha
  \text{ is represented by }
  \sum_{i}
  \lambda_i
  \sigma_i
  \right\}
\]
For a closed connected and oriented manifold $\Manifold$ the norm of the
fundamental class
$\Norm{\FundamentalClass{\Manifold}}$ is called the simplicial volume of
$\Manifold$ and denoted by $\SimplicialVolume{\Manifold}$. This was first
investigated by Gromov in \cite{Gromov} and has various relations to topology
and geometry. It is a consequence of Gromovs Mapping Theorem
\cite[The Mapping Theorem, Section 3.1]{Gromov} that the simplicial volume of a
manifold only
depends on the image of the fundamental class in the group homology of its
fundamental group with respect to the map that classifies the universal
covering.

Another consequence of Gromovs Mapping Theorem is that the simplicial
volume of every manifold $\ManifoldTotal$, which fibers over a manifold
$\ManifoldBase$ with a fiber that has amenable fundamental group and has
positive dimension, vanishes. Interestingly, if one instead
requires that the base manifold has amenable fundamental group, while having no
requirements on the fiber, then an analogue of this
result can not hold. In
fact every mapping torus of a pseudo-Anosov surface automorphism is a
hyperbolic $3$-manifold and hence has non-zero simplicial volume. Strengthening
the restrictions on the base manifold leads to
the following question:
\begin{question}
\label{qst:Base1-con}
   	Suppose a closed, connected and oriented manifold $\ManifoldTotal$ fibers
   	over a simply-connected manifold $\ManifoldBase$ with positive dimension,
   	does the simplicial volume of $\ManifoldTotal$ vanish?
\end{question}

As a partial answer to \autoref{qst:Base1-con} we are able to show the
following:
\input{Corollary/Vanishing.tex}
Note that if the base manifold is simply-connected, then the structure group of
the bundle can always be reduced to a connected group (See
\autoref{scn:FBandSG} for an introduction to structure groups). Hence in the
context of
\autoref{qst:Base1-con} it is quite natural to consider bundles with connected
structure group and no restriction on the fundamental group of the base
manifold as in \autoref{cor:MainCorollary}.

As a previous result, it was shown in \cite{KastenholzReinhold} that if
$\ManifoldBase$ is either a
simply-connected sphere, or a flexible and $2$-connected manifold and the
fundamental group of the fiber satisfies the Novikov conjecture, then the
simplicial volume of the total space vanishes.

The structure group being compact is equivalent to the existence of a
metric on the fiber which is invariant under the structure group. It is futher
equivalent to the existence of a
metric on the total space such that the fibers are totally geodesic
submanifolds (See for example Theorem 1.11 and Theorem 1.12 in
\cite{TotallyGeodesic}). Furthermore since isometry groups are
Lie groups, this shows that the first case is actually contained in the second
case (if one passes to an effective action of the structure group).
Additionally every Lie group is homotopy equivalent to its maximally
compact structure group, hence the first and second point are actually
equivalent.

\autoref{cor:MainCorollary} will be a consequence of the following more general
result:
\input{Theorem/RationalProductVanishing.tex}

The proof of \autoref{thm:RationalProductVanishing} can roughly be divided into
two parts. Let $\FundamentalGroupShortHand$ denote the fundamental group of the
base manifold. In the first part we show that the simplicial volume of the
total space only depends on the image of the fundamental class of the base
manifold under the canonical map to
$
   	\ClassifyingSpace{\Group}
   	\times
   	\ClassifyingSpace{\FundamentalGroupShortHand}
$
(See \autoref{thm:HomologyInvariance} for the precise statement). This step
only requires that $\Group$ is connected.

For the second step, we will establish in \autoref{scn:Representatives} that
homology classes in $\ClassifyingSpace{\Group}$ can be represented by a product
of simply-connected spheres, provided that $\Group$ is connected and
$\ClassifyingSpace{\Group}$ has the rational homotopy type of a product of
Eilenberg-MacLane spaces and then combine these representatives with Theorem~C
in \cite{KastenholzReinhold} to finish the proof of
\autoref{thm:RationalProductVanishing}. Finally we will also show in
\autoref{scn:Representatives} that the classifying spaces of the structure
groups in \autoref{cor:MainCorollary} have the rational homotopy type of a
product of Eilenberg-MacLane spaces.

\begin{remark}
  The proof presented in \autoref{scn:Proof} actually shows that if the answer
  to \autoref{qst:Base1-con} is yes, then the simplicial volume of the total
  space of any bundle with connected structure group agrees with the simplicial
  volume of the trivial bundle.
\end{remark}

\subsection*{Acknowledgements}
	The author would like to thank Mark Pedron and Manuel Krannich for helpful
	discussions about rational homotopy theory.

%% file: Corollary/Vanishing.tex
\begin{theorem}
\label{cor:MainCorollary}
  Let $\ManifoldFiber \to \ManifoldTotal \to \ManifoldBase$ denote a bundle of closed, connected and oriented manifolds with either
  \begin{enumerate}[(a)]
   \item
   	compact structure group, or
   \item 
   	a connected Lie group as structure group, or
   \item
   	 aspherical fiber and connected structure group
  \end{enumerate} 	
  then the simplicial volume of
  $\ManifoldTotal$ agrees with the simplicial volume of the trivial
  bundle $\ManifoldFiber \times \ManifoldBase$.
\end{theorem}

%% file: Theorem/RationalProductVanishing.tex
\begin{theorem}
\label{thm:RationalProductVanishing}
  Let $\ManifoldFiber \to \ManifoldTotal \to \ManifoldBase$ denote a fiber
  bundle of closed oriented connected manifolds with structure group $\Group$.
  If $\HomotopyGroupOfObject{0}{\Group}{e}$ is finite and a multiple of every 
  rational homology class in the universal covering 
  $\UniversalCovering{\ClassifyingSpace{\Group}}$ can be 
  represented by a disjoint union of spaces that fiber over simply-connected 
  spheres, then the simplicial volume of $\ManifoldTotal$
  agrees with the simplicial volume of the trivial bundle $\ManifoldBase \times
  \ManifoldFiber$.
  
  In particular this holds if $\UniversalCovering{\ClassifyingSpace{\Group}}$ 
  has the rational homotopy type of a product of simply-connected
  Eilenberg-Maclane spaces, 
\end{theorem}

%% file: Section/FBandSG.tex
\label{scn:FBandSG}
This section contains a short introduction to fiber bundles and their structure
groups, for a more thorough introduction see \cite{SomethingFiberBundle}.

Let $\FiberingProjektion{\ManifoldTotal}\colon \ManifoldTotal \to
\ManifoldBase$ denote a smooth fiber
bundle with fiber $\ManifoldFiber$, i.e. there exists an open covering $U_i
\subset \ManifoldBase$ and fiber preserving
diffeomorphisms
$
  \Phi_i
  \colon
  \apply
    {\FiberingProjektion{\ManifoldFiber}^{-1}}
    {U_i}
  \to
  U_i
  \times
  \ManifoldFiber
$
such that
$
  \Phi_i \circ \Phi_j^{-1}
  \colon
  U_i \cap U_j
  \times
  \ManifoldFiber
  \to
  U_i \cap U_j
  \times
  \ManifoldFiber
$
is a fiber preserving diffeomorphism.
In particular this yields maps
$
  \overline{\Phi_{i,j}}
  \colon
  U_i \cap U_j
  \to
  \Diff{\ManifoldFiber}
$
such that
\begin{equation}
\label{eqn:CocycleCondition}
  \overline{\Phi_{k,i}} \circ \overline{\Phi_{i,j}}
  =
  \overline{\Phi_{k,j}}
\end{equation}%

These maps are usually called \introduce{cocycles of $\ManifoldTotal$} and
since these cocycles carry the gluing information of the bundle, one can check
that there is a one-to-one correspondence between maps satisfying the cocycle
condition (\ref{eqn:CocycleCondition}) and fiber bundles.
Notice that the same cocycles can be used to define a
$\Diff{\ManifoldFiber}$-principal bundle
$\ManifoldTotal_{\Diff{\ManifoldFiber}}$ which comes equipped with a map
$
  \ManifoldTotal_{\Diff{\ManifoldFiber}}
  \times
  \ManifoldFiber
  \to
  \ManifoldTotal
$%
. Conversely, using the left action of $\Diff{\ManifoldFiber}$ on
$\ManifoldFiber$, one can assosciate to every $\Diff{\ManifoldFiber}$-principal
bundle $P$ a $\ManifoldFiber$-bundle by taking the fiber product $P \times
\ManifoldFiber/ \Diff{\ManifoldFiber}$.
This yields a one-to-one correspondence between $\Diff{\ManifoldFiber}$-bundles
and $\ManifoldFiber$-bundles.

There is the universal $\Diff{\ManifoldFiber}$-principal bundle
$
  \Diff{\ManifoldFiber}
  \to
  E\Diff{\ManifoldFiber}
  \to
  \ClassifyingSpace{\Diff{\ManifoldFiber}}
$
and for every $\Diff{\ManifoldFiber}$-principal bundle $P$ over
$\ManifoldBase$, there exists a up to homotopy unique map
$
  \ClassifyingMap{P}
  \colon
  \ManifoldBase
  \to
  \ClassifyingSpace{\Diff{\ManifoldFiber}}
$
such that the pullback $\ClassifyingMap{P}^* E\Diff{\ManifoldFiber}$ is
isomorphic to $P$. Analogously, for every $\ManifoldFiber$-bundle
$\ManifoldTotal$ over $\ManifoldBase$, there exists a up to homotopy unique map
$
  \ClassifyingMap{\ManifoldTotal}
  \colon
  \ManifoldBase
  \to
  \ClassifyingSpace{\Diff{\ManifoldFiber}}
$
such that $\ManifoldTotal$ is isomorphic to
$\ClassifyingMap{\ManifoldTotal}^*\TautologicalBundle{\ManifoldFiber}$, with
$
  \TautologicalBundle{\ManifoldFiber}
  =
  E\Diff{\ManifoldFiber} \times_{\Diff{\ManifoldFiber}} \ManifoldFiber
$%
. For that reason $\ClassifyingSpace{\Diff{\ManifoldFiber}}$ is called the
classifying space for $\Diff{\ManifoldFiber}$ and we call
$\TautologicalBundle{\ManifoldFiber}$, the universal $\ManifoldFiber$ bundle.

Let now $\Group$ denote a topological group admitting a continuous group action
on $\ManifoldFiber$.
Such a group action is equivalent to a continuous group homomorphism
\[
    \phi \colon \Group \to \Diff{\ManifoldFiber}
\]
We say that \introduce{the structure group of a bundle $\ManifoldTotal \to
\ManifoldBase$
can be reduced to $\Group$}, if there exists an open covering of $\ManifoldBase$
and trivializations over these open maps, such that the cocycles
$
  \overline{\Phi_{i,j}}\colon U_i \cap U_j \to \Diff{\ManifoldFiber}
$
can be lifted to maps $\widetilde{\Phi_{i,j}} \colon U_i \cap U_j \to \Group$
such that the coycle condition
$
  \widetilde{\Phi_{k,i}} \circ \widetilde{\Phi_{i,j}}
  =
  \widetilde{\Phi_{k,j}}
$%
still holds. Am $\ManifoldFiber$-bundle together with a choice of these lifts is
called a \introduce{$\ManifoldFiber$-bundle with structure group $\Group$}.

The group homomorphism $\phi$ also yields a map on classifying spaces
$
    \ClassifyingSpace{\phi}
    \colon
    \ClassifyingSpace{\Group}
    \to
    \ClassifyingSpace{\Diff{\ManifoldFiber}}
$
and the structure group of a bundle can be reduced to $\Group$ if and only if
the classifying map lifts to $\ClassifyingSpace{\Group}$ and a reduction of the
structure group corresponds to such a lift. Similar to the previous
considerations, for every $\ManifoldFiber$-bundle $\ManifoldTotal \to
\ManifoldBase$ with structure group $\Group$, there exists a up to homotopy
unique map $\ClassifyingMap{\ManifoldTotal}\colon \ManifoldBase \to
\ClassifyingSpace{\Group}$ such that $\ManifoldTotal$ is isomorphic as a bundle
with structure group $\Group$ (meaning that the fiberwise isomorphisms are also
elements in $\Group$) to
$
  \ClassifyingMap{\ManifoldTotal}^*
  \left(
    \ClassifyingSpace{\phi}^* \TautologicalBundle{\ManifoldFiber}
  \right)
$%
. We denote $\ClassifyingSpace{\phi}^* \TautologicalBundle{\ManifoldFiber}$ by
$\TautologicalBundle{\ManifoldFiber,\Group}$ and call it the universal
$\ManifoldFiber$-bundle with structure group $\Group$.

\begin{example}
  Let $\ManifoldTotal \to \ManifoldBase$ denote a smooth $\Reals^n$, bundle,
  then a reduction of the structure group of $\ManifoldTotal$ to
  $\GeneralLinearGroup{n} \subset \Diff{\Reals^n}$ equips $\ManifoldTotal$ with
  a vector bundle structure. Furthermore the structure group of such a bundle
  can be reduced to $\apply{\mathop{GL}^+_{n}}{\Reals}$ if and only if the
  bundle is orientable and a choice of orientation corresponds to a reduction
  of the structure group to $\apply{\mathop{GL}^+_{n}}{\Reals}$.

  One can proceed similarly with spin-structures, or metrics, or
  complex structures.
\end{example}

Another important example for us is the following:
\begin{example}
  Let $\Group$ denote a topological group and $\Group_0$ its identity
  component, then one easily sees that
  $
    \ClassifyingSpace{\Group_0}
    \to
    \ClassifyingSpace{\Group}
  $
  is a universal covering map. As a consequence, one can reduce the structure
  group of every bundle $\ManifoldTotal \to \ManifoldBase$ with structure group
  $\Group$ and simply-connected base $\ManifoldBase$ to $\Group_0$.
\end{example}

Note that the very same considerations will also work for topological bundles
i.e. bundles with structure group $\HomeoGroup{\ManifoldFiber}$ and even
fibrations, where one takes the assosciative H-space of self-homotopy
equivalences $\Automorphisms{\ManifoldFiber}$ as $\Group$ (and the cocycle
condition only has to hold up to homotopy).

%% file: Section/Proof.tex
\label{scn:Proof}
The following section contains the proof of
\autoref{thm:RationalProductVanishing}. The proof consists of two parts, first
we show a connection between the simplicial volume of the total space of a
bundle and the homology classes in
$
	\ClassifyingSpace{\Group}
	\times
	\ClassifyingSpace{\FundamentalGroupShortHand}
$%
. This step works for every connected structure group.
Finally we use the existence of simple representatives of all homology class
in $\ClassifyingSpace{\Group}$ to finish the proof.

Let $\ManifoldFiber \to \ManifoldTotal\to \ManifoldBase$ denote a fiber bundle
of compact, connected and orientable manifolds
with connected structure group $\Group$ and $n$-dimensional base
manifold.

Let
$
  \ClassifyingMap{\ManifoldTotal}
  \colon
  \ManifoldBase
  \to
  \ClassifyingSpace{\Group}
$
denote the classifying map of the bundle, i.e. there exists an isomorphism
$
  \ClassifyingMap{\ManifoldTotal}^* \TautologicalBundle{\ManifoldFiber,\Group}
  \cong
  \ManifoldTotal
$%
, where $\TautologicalBundle{\ManifoldFiber,\Group}$ denotes the tautological
$\ManifoldFiber$ bundle with structure group $\Group$ over
$\ClassifyingSpace{\Group}$
as introduced in \autoref{scn:FBandSG}.

Let us denote the fundamental
group of $\ManifoldBase$ by $\FundamentalGroupShortHand$, and let
$\ClassifyingMap{\UniversalCovering{\ManifoldBase}}\colon \ManifoldBase \to
\ClassifyingSpace{\FundamentalGroupShortHand}$ denote a classifying map of the
universal covering. Consider the map
\[
  \ContinuousMap
  =
  \ClassifyingMap{\ManifoldTotal}
  \times
  \ClassifyingMap{\UniversalCovering{\ManifoldBase}}
  \colon
  \ManifoldBase
  \to
  \ClassifyingSpace{\Group}
  \times
  \ClassifyingSpace{\FundamentalGroupShortHand}
\]
which is covered by a map
$
  \overline{\ContinuousMap}
  \colon
  \ManifoldTotal
  \to
  \TautologicalBundle{\ManifoldFiber,\Group}
  \times
  \ClassifyingSpace{\FundamentalGroupShortHand}
$%
. We denote the projection
$
  \TautologicalBundle{\ManifoldFiber,\Group}
  \times
  \ClassifyingSpace{\FundamentalGroupShortHand}
  \to
  \ClassifyingSpace{\Group}
  \times
  \ClassifyingSpace{\FundamentalGroupShortHand}
$
by $\FiberingProjektion{\TautologicalBundle{\ManifoldFiber,\Group}}$.

The corresponding long exact sequences of the homotopy groups of a fibration
yield the following diagram:
\begin{center}
  \begin{tikzcd}
    \ldots
      \ar[r]
    &
    \HomotopyGroupOfObject{2}{\ManifoldBase}{\Point}
      \ar[r]
      \ar[d]
    &\ClassifyingMap{\UniversalCovering{\ManifoldBase}}
    \HomotopyGroupOfObject{1}{\ManifoldFiber}{\Point}
      \ar[r]
      \ar[d,"="]
    &
    \HomotopyGroupOfObject{1}{\ManifoldTotal}{\Point}
      \ar[r]
      \ar[d,"\HomotopyGroupMorphism{\overline{\ContinuousMap}}"]
    &
    \HomotopyGroupOfObject{1}{\ManifoldBase}{\Point}
      \ar[r]
      \ar[d,"\cong"]
    &
    0
    \\
    \ldots
      \ar[r]
    &
    \HomotopyGroupOfObject
      {2}
      {
        \ClassifyingSpace{\Group}
        \times
        \ClassifyingSpace{\FundamentalGroupShortHand}
      }
      {\Point}
      \ar[r]
    &
    \HomotopyGroupOfObject
      {1}
      {\ManifoldFiber}
      {\Point}
      \ar[r]
    &
    \HomotopyGroupOfObject
      {1}
      {
        \TautologicalBundle{\ManifoldFiber,\Group}
        \times
        \ClassifyingSpace{\FundamentalGroupShortHand}
      }
      {\Point}
      \ar[r]
    &
    \HomotopyGroupOfObject
      {1}
      {
        \ClassifyingSpace{\Group}
        \times
        \ClassifyingSpace{\FundamentalGroupShortHand}
      }
      {\Point}
      \ar[r]
    &
    0
  \end{tikzcd}
\end{center}
A simple diagram chase shows that
$\HomotopyGroupMorphism{\overline{\ContinuousMap}}$
is a surjection with amenable kernel, in particular by Gromovs mapping theorem,
$\overline{\ContinuousMap}$ induces an isometry on homology with respect to the
$\ellone$-norm.
Therefore the simplicial volume of $\ManifoldTotal$ agrees with the
$\ellone$-norm of
$
  \apply{\overline{\ContinuousMap}}{\FundamentalClass{\ManifoldTotal}}
$,
which equals
$
  \apply{
    \FiberTransfer{\FiberingProjektion{\TautologicalBundle{\ManifoldFiber,\Group}}}
    }
    {
      \apply
        {\HomologyOfSpaceMorphism{\ContinuousMap}}
        {\FundamentalClass{\ManifoldBase}}
    }
$%
, where
$\FiberTransfer{\FiberingProjektion{\TautologicalBundle{\ManifoldFiber,\Group}}}$
denotes the fiber transfer. In particular
$
  \apply{
    \FiberTransfer{\FiberingProjektion{\TautologicalBundle{\ManifoldFiber,\Group}}}
  }
  {
    \apply
    {\HomologyOfSpaceMorphism{\ContinuousMap}}
    {\FundamentalClass{\ManifoldBase}}
  }
$
only depends on $\apply{\ContinuousMap}{\FundamentalClass{\ManifoldBase}}$.
All in all this yields
\input{Theorem/HomologyInvariance.tex}
We will now focus on classifying space admitting nice representatives of their
homology classes. This is encapsulated in the following definition:
\input{Definition/RationallySpherical.tex}

Let us from here on forth assume that the rational homology of
$\ClassifyingSpace{\Group}$ is generated by fiber bundles over spheres. We
will show in
\autoref{scn:Representatives} that this includes the three
cases:
\begin{itemize}
\item
	$\Group$ is a compact group
\item
  $\Group$ is a Lie group
\item
  $\Group$ is the identity component of the H-space of self-homotopy
  equivalences of an aspherical manifold.
\end{itemize}

By Künneths theorem, the $n$-th rational homology group of
$
  \ClassifyingSpace{\Group}
  \times
  \ClassifyingSpace{\FundamentalGroupShortHand}
$
is isomorphic to
\[
  \bigoplus_{i}
  \HomologyOfGroupObject{i}{\ClassifyingSpace{\Group}}{\Rationals}
  \otimes
  \HomologyOfGroupObject
    {n-i}
    {\ClassifyingSpace{\FundamentalGroupShortHand}}
    {\Rationals}
\]
Hence
$
  \apply
    {\HomologyOfSpaceMorphism{\ContinuousMap}}
    {\FundamentalClass{\ManifoldBase}}
$
is the sum of classes $\alpha_i$ with
\[
  \alpha_i
  \in
  \HomologyOfGroupObject{i}{\ClassifyingSpace{\Group}}{\Rationals}
  \otimes
  \HomologyOfGroupObject
    {n-i}
    {\ClassifyingSpace{\FundamentalGroupShortHand}}
    {\Rationals}
\]

By assumption any $\alpha_i$ with $i\geq 1$
can be represented by a disjoint union of manifolds fibering over
simply-connected spheres times some representative of a class in
$
  \HomologyOfGroupObject
    {n-i}
    {\ClassifyingSpace{\FundamentalGroupShortHand}}
    {\Rationals}
$%
. Therefore the norm of
$
  \apply{
    \FiberTransfer{\FiberingProjektion{\TautologicalBundle{\ManifoldFiber,\Group}}}
  }
  {\alpha_i}
$
which equals the simplicial volume of the pullback of
$\TautologicalBundle{\ManifoldFiber,\Group}$ to the aforementioned
representative
vanishes by Theorem~C in \cite{KastenholzReinhold} as the simplicial
volume of any space fibering over a simply-connected sphere vanishes. Therefore
we conclude that the simplicial volume of $\ManifoldTotal$ agrees with the
$\ellone$-norm of
$
  \apply{
      \FiberTransfer{\FiberingProjektion{\TautologicalBundle{\ManifoldFiber,\Group}}}
    }
    {\alpha_0}
$%
, but $\alpha_0$ is an element of
$
    \HomologyOfGroupObject{0}{\ClassifyingSpace{\Group}}{\Rationals}
    \otimes
    \HomologyOfGroupObject
      {n}
      {\ClassifyingSpace{\FundamentalGroupShortHand}}
      {\Rationals}
$
to be more precise it equals
$
  \apply
    {
      \HomologyOfSpaceMorphism{\ClassifyingMap{\UniversalCovering{\ManifoldBase}}}
    }
    {\FundamentalClass{\ManifoldBase}}
$
, in particular it is represented by
$
  \ast
  \times
  \ClassifyingMap{\UniversalCovering{\ManifoldBase}}
  =
  \ClassifyingMap{\ManifoldFiber \times \ManifoldBase}
  \times
  \ClassifyingMap{\UniversalCovering{\ManifoldBase}}
  \colon
  \ManifoldBase
  \to
  \ClassifyingSpace{\Group}
  \times
  \ClassifyingSpace{\FundamentalGroupShortHand}
$
and therefore
$
  \apply
    {
      \FiberTransfer{
        \FiberingProjektion{\TautologicalBundle{\ManifoldFiber,\Group}}
      }
    }
    {\alpha_0}
$
is represented by $\ManifoldFiber \times \ManifoldBase$ which concludes the
proof of \autoref{thm:RationalProductVanishing}.

%% file: Theorem/HomologyInvariance.tex
\begin{proposition}
\label{thm:HomologyInvariance}
  Let
  $
    \ManifoldFiber
    \to
    \ManifoldTotal
    \to
    \ManifoldBase
  $
  and
  $
    \ManifoldFiber
    \to
    \ManifoldTotal'
    \to
    \ManifoldBase'
  $
  denote two bundles with connected structure group $\Group$ and suppose
  further that the fundamental groups of $\ManifoldBase$ and $\ManifoldBase'$
  are isomorphic. Let us denote their fundamental group by
  $\FundamentalGroupShortHand$.
  Assume further that there exists a $\lambda \in \Reals$ such that
  \[
    \apply{
      \left(
        \ClassifyingMap{\ManifoldTotal}
        \times
        \ClassifyingMap{\UniversalCovering{\ManifoldBase}}
      \right)_*
    }
    {\FundamentalClass{\ManifoldBase}}
    =
    \lambda
    \apply{
      \left(
        \ClassifyingMap{\ManifoldTotal'}
        \times
        \ClassifyingMap{\UniversalCovering{\ManifoldBase'}}
      \right)_*
      }
    {\FundamentalClass{\ManifoldBase'}}
    \in
    \HomologyOfSpaceObject
      {*}
      {
        \ClassifyingSpace{\Group}
        \times
        \ClassifyingSpace{\FundamentalGroupShortHand}
      }
      {\Reals}
  \]
  then
  $
    \SimplicialVolume{\ManifoldTotal}
    =
    \lambda \SimplicialVolume{\ManifoldTotal'}
  $.
\end{proposition}

%% file: Definition/RationallySpherical.tex
\begin{definition}
	Let $\TopologicalSpace$ denote a simply-connected topological space, we
	call its rational homology \introduce{generated by fiber bundles over
	spheres} if
	there exists maps
	\[
		\ContinuousMap_i
		\colon
    \ManifoldTotal_i
		\to
		\TopologicalSpace
	\]
	where $i$ lies in some index set and $\ManifoldTotal_i$ is a closed, oriented
	manifold that fibers over a simply-connected sphere, such that the set
	$
		\left\{
			\apply
                {\HomologyOfSpaceMorphism{\ContinuousMap_{i}}}
                {
                    \FundamentalClass
                        {\ManifoldTotal_i}
                }
		\right\}
	$
    generates the rational homology groups of $\TopologicalSpace$.
\end{definition}

%% file: Section/Representatives.tex
\label{scn:Representatives}
The goal of this section is to show that the rational homology of the
classifying space $\ClassifyingSpace{\Group}$ is generated by fiber bundles
over spheres if:
\begin{itemize}
  \item
  $\Group$ is a compact group
  \item
  $\Group$ is a Lie group
  \item
  $\Group$ is the identity component of the H-space of self-homotopy
  equivalences of an aspherical manifold.
\end{itemize}

In order to accomplish that, we will show that the rational
homology of spaces that are
rationally homotopy equivalent to a product of Eilenberg-MacLane spaces is
generated by fiber bundles over spheres. In fact, we will show something
slightly stronger:
\input{Definition/SufficientlyStronglySpherical.tex}

The following lemma shows that this property is preserved under rational
homotopy equivalences. This follows immediately from the fact that products of
spheres are universal (See \cite{UniversalSpaces}).
Nevertheless we include a proof here.
\input{Lemma/LiftabilityRationalEquivalence.tex}

Let us now restrict to spaces rationally equivalent to products of rational
Eilenberg-Maclane spaces for which
it is shown quite easily that their rational homology is generated by products
of spheres.
\input{Proposition/RepresentativesProductEMSpaces.tex}

The following proposition ties the previous considerations to the main cases in
question.
\input{Proposition/ProductEMSpaces.tex}
As a consequence, we immediately obtain the following, which combined with
\autoref{thm:RationalProductVanishing} yields \autoref{cor:MainCorollary}.
\input{Corollary/RepresentativesBG.tex}

%% file: Definition/SufficientlyStronglySpherical.tex
\begin{definition}
  Let $\TopologicalSpace$ denote a simply-connected topological space, we
  call its rational homology \introduce{generated by products of spheres} if
  there exists maps
  \[
    \ContinuousMap_i
    \colon
    \prod_{l=1}^{k_{i}}
    \Sphere{n_{i,l}}
    \to
    \TopologicalSpace
  \]
  where $i$ lies in some index set, and $n_{i,l}\geq 2$, such that
  $
    \left\{
    \apply
    {\HomologyOfSpaceMorphism{\ContinuousMap_{i}}}
    {
      \FundamentalClass
      {\prod_{l=1}^{k_{i}} \Sphere{n_{i,l}}}
    }
    \right\}
  $
  generates the rational homology groups of $\TopologicalSpace$.
\end{definition}

%% file: Lemma/LiftabilityRationalEquivalence.tex
\begin{lemma}
\label{lem:Liftability}
  Suppose $\ContinuousMap \colon \TopologicalSpace \to \TopologicalSpace'$ is a
  rational homotopy equivalence and suppose further that the rational homology
  of $\TopologicalSpace'$ is generated by products of spheres, then the same
  holds for $\TopologicalSpace$.
\end{lemma}
\begin{proof}
  Let $S = \prod_{i} \Sphere{n_i}$ denote a product of spheres and
  $\ContinuousMapALT \colon S \to \TopologicalSpace'$ a map. We will show that
  there exists a map $h\colon S \to S$ of non-zero degree such that
  $\ContinuousMapALT \circ h$ lifts along $\ContinuousMap$.

  Note that $\ContinuousMap$ being a rational homotopy equivalence implies that
  all homotopy groups of the homotopy fiber $\HomotopyFiber{\ContinuousMap}$
  are torsion. We will construct a lift of $\ContinuousMapALT$ via induction
  over its skeleta $S^{(d)}$, precomposing with a map of non-zero degree map
  where necessary.

  Let $h_k\colon S \to S$ denote the map that is the product of
  degree $k$ maps on the factors of $S$. By the Künneth theorem, this map sends
  every integral homology class to a $k$-fold multiple of another integral
  homology class, additionally we have $h_k \circ h_{k'} = h_{kk'}$.

  Since $\ContinuousMap$ induces an isomorphism on rational homotopy groups,
  there exists a $k_0$ such that, we
  can lift $\at{\ContinuousMapALT\circ h_{k_0}}{\bigvee_i \Sphere{n_i}}$ to
  $\TopologicalSpace$.
  In particular we can lift $\ContinuousMapALT\circ h_{k_0}$ to the first
  $d$-skeleton which is not a point.

  Let us now proceed inductively and assume that a lift on the $d$ skeleton of
  $S$ of $\ContinuousMapALT \circ h_{k_d}$ was already
  constructed.
  The obstruction class
  $o_{d+1}$ to extending this lift to the $d+1$ skeleton is an element of
  \[
    \CohomologyOfSpaceObject
      {d+1}
      {S^{(d+1)},S^{(d)}}
      {\HomotopyGroupOfObject{d}{\HomotopyFiber{\ContinuousMap}}{\Point}}
  \]
  Let us consider the long exact sequence of integral cohomology for this pair:
  \[
    0
    \to
    \CohomologyOfSpaceObject{d}{S^{(d+1)}}{\Integers}
    \to
    \CohomologyOfSpaceObject{d}{S^{(d)}}{\Integers}
    \to
    \CohomologyOfSpaceObject{d+1}{S^{(d+1)},S^{(d)}}{\Integers}
    \to
    \CohomologyOfSpaceObject{d+1}{S^{(d+1)}}{\Integers}
    \to
    0
  \]
  Since the cellular chain complex of a sphere is concentrated in a single
  degree, one concludes that the map
  $
    \CohomologyOfSpaceObject{d}{S^{(d+1)}}{\Integers}
    \to
    \CohomologyOfSpaceObject{d}{S^{(d)}}{\Integers}
  $
  is an isomorphism, hence
  $\CohomologyOfSpaceObject{d+1}{S^{(d+1)},S^{(d)}}{\Integers}$ is isomorphic
  to
  $\CohomologyOfSpaceObject{d+1}{S^{(d+1)}}{\Integers}$. Therefore $h_k$ also
  sends every class in $\CohomologyOfSpaceObject{d+1}{S^{(d+1)}}{\Integers}$ to
  a $k$-fold multiple of another class. Hence there exists a $k_{d+1}$ such
  that
  $
    \apply{h_{k_{d+1}}^{*}}{o_{d+1}}
  $
  vanishes. By naturality of the obstruction class we conclude that the lift of
  $
    \ContinuousMapALT\circ h_{k_d} \circ h_{k_{d+1}}
    =
    \ContinuousMapALT \circ h_{k_d k_{d+1}}
  $
  can be extended to the $d+1$-skeleton $S^{(d+1)}$.
\end{proof}

%% file: Proposition/RepresentativesProductEMSpaces.tex
\begin{proposition}
\label{prp:RepresentativesRationalHomotoptyProdEMSpaces}
  Let $\TopologicalSpace$ denote a space that has the rational homotopy type of
  a product of simply-connected rational Eilenberg-Maclane
  spaces, then its rational homology is generated by products of spheres.
\end{proposition}
\begin{proof}
  Note that the rational cohomology ring of $\EMSpace{\Rationals}{n}$ is either
  an exterior algebra with a single generator in degree $n$ if $n$ is odd, or a
  polynomial algebra generated by a single generator in degree $n$ if $n$ is
  even. Evidently, if $n$ is odd, then every rational homology class is
  representable by spheres. If $n$ is even, consider a generator in degree $kn$
  and let $\iota_k\colon \left(\Sphere{n}\right)^{k} \to
  \EMSpace{\Rationals}{n}$,
  denote the map that classifies
  $
    \alpha
    =
    \sum_{i=1}^{k}
    \apply{\pi_i^*}{\FundamentalClass{\Sphere{n}}^*}
  $
  where $\pi_i$ is the projection onto the $i$-th factor in
  $\left(\Sphere{n}\right)^{k}$ and $\FundamentalClass{\Sphere{n}}^*$ denotes
  the dual fundamental class. Since $\alpha^k$ is a non-zero multiple of the
  dual fundamental class of $\left(\Sphere{n}\right)^k$, it follows that
  $\iota_k$ induces an isomorphism in rational cohomology in degree $kn$.
  Hence
  $\iota_k$ represents a multiple of the generator in the $kn$-th homology of
  $\EMSpace{\Rationals}{n}$.

  Since the rational homology groups of a product are generated by cross
  products of classes in the factors it follows that any rational homology
  class in the product of simply-connected rational Eilenberg-Maclane spaces
  can be represented as a sum of classes representable by products of spheres.

  Therefore \autoref{lem:Liftability} yields the desired result.
\end{proof}

%% file: Proposition/ProductEMSpaces.tex
\begin{proposition}
\label{prp:ProductEMSpaces}
  Let $\CompactGroup$ denote a connected topological group whose underlying
  topological space is of finite type, then the classifying
  space $\ClassifyingSpace{\CompactGroup}$ has the rational homotopy type of a
  product of simply-connected rational Eilenberg-MacLane spaces. This class
  evidently contains compact groups and Lie groups.

  If $\TopologicalSpace$ is an aspherical topological space, then the
  classifying space of the identity component of its space of self-homotopy
  equivalences
  $\ClassifyingSpace{\AutomorphismsConnected{\TopologicalSpace}}$ is rationally
  homotopy equivalent to a $\EMSpace{V}{2}$, where
  $
    V
    =
    \apply
      {Z}
      {\HomotopyGroupOfObject{1}{\TopologicalSpace}{\Point}}
    \otimes
    \Rationals
  $%
  and
  $
    \apply
      {Z}
      {\HomotopyGroupOfObject{1}{\TopologicalSpace}{\Point}}
  $ denotes the center of the fundamental group of $\TopologicalSpace$.

\end{proposition}
\begin{proof}
  For the first part, note that by Borel-Serres theorem (See for example
  \cite[Theorem 6.38]{McCleary}) the cohomology ring of
  $\ClassifyingSpace{\CompactGroup}$ is a polynomial ring with finitely many
  generators in even degree. Therefore the product of the maps classifying the
  generators
  $
    \ClassifyingSpace{\CompactGroup}
    \to
    \prod_{i}
    \EMSpace{\Rationals}{n_{i}}
  $
  induces an isomorphism in rational cohomology by the Künneth theorem and
  hence a rational homotopy equivalence.

  For the second part note that
  $\HomotopyGroupOfObject{i}{\Automorphisms{\TopologicalSpace}}{\Identity}$ is
  $
    \apply
      {\mathop{Out}}
      {\HomotopyGroupOfObject{1}{\TopologicalSpace}{\Point}}
  $
  if $i=0$,
  $
    \apply
      {Z}
      {\HomotopyGroupOfObject{1}{\TopologicalSpace}{\Point}}
  $
  if $i=1$ and zero else. Therefore
  $\ClassifyingSpace{\AutomorphismsConnected{\TopologicalSpace}}$ is already a
  $\EMSpace
    {
      \apply
        {Z}
        {\HomotopyGroupOfObject{1}{\TopologicalSpace}{\Point}}
    }
    {2}
  $%
  .
\end{proof}

%% file: Corollary/RepresentativesBG.tex
\begin{corollary}
\label{cor:RepresentantivesBG}
  Let $\Group$ denote either a connected compact group, a Lie group or the
  identity component of the homotopy automorphisms of an aspherical topological
  space, then the rational homology of $\ClassifyingSpace{\Group}$ is
  generated by products of spheres, in particular it is generated by
  fiber bundles over spheres.
\end{corollary}

%% file: top.bbl
\begin{thebibliography}{1}

\bibitem{TotallyGeodesic}
Richard H.~jun. Escobales.
\newblock Riemannian submersions with totally geodesic fibers.
\newblock {\em J. Differ. Geom.}, 10:253--276, 1975.

\bibitem{Gromov}
M.~Gromov.
\newblock Volume and bounded cohomology.
\newblock {\em Publ. Math., Inst. Hautes {\'E}tud. Sci.}, 56:5--99, 1982.

\bibitem{SomethingFiberBundle}
Dale~H. Husemoller.
\newblock {\em Fibre bundles.}, volume~20 of {\em Grad. Texts Math.}
\newblock Berlin: Springer-Verlag, 3rd ed. edition, 1993.

\bibitem{KastenholzReinhold}
Thorben Kastenholz and Jens Reinhold.
\newblock Simplicial volume and essentiality of manifolds fibered over spheres.
\newblock {\em J. Topol.}, 16(1):192--206, 2023.

\bibitem{McCleary}
John McCleary.
\newblock {\em A user's guide to spectral sequences.}, volume~58 of {\em Camb.
  Stud. Adv. Math.}
\newblock Cambridge: Cambridge University Press, 2nd ed. edition, 2001.

\bibitem{UniversalSpaces}
Mamoru Mimura, Ronald~C. O'Neill, and Hirosi Toda.
\newblock On p-equivalence in the sense of {Serre}.
\newblock {\em Jpn. J. Math.}, 40:1--10, 1971.

\end{thebibliography}
